\documentclass[12pt,a4paper,leqno]{article}

\usepackage{amssymb,exscale}
\usepackage[centertags]{amsmath}
\usepackage{graphicx}
\usepackage{theorem}

\numberwithin{equation}{section}

\theoremstyle{change}
\theorembodyfont{\rmfamily}

\newtheorem{lemma}[equation]{Lemma}
\newtheorem{proposition}[equation]{Proposition}

\newtheorem{problem}[equation]{Problem}

\newcommand{\sk}{\smallskip}

\renewcommand{\phi}{\varphi}

\renewcommand{\(}{\bigl(}
\renewcommand{\)}{\bigr)\vphantom{)}}

\newcommand{\Om}{\Omega}

\newcommand{\de}{\delta}

\newcommand{\E}{\mathbb E}
\newcommand{\Ex}[1]{\,\mathbb E\,\(\,#1\,\)\,}
\newcommand{\Pro}[1]{\,\mathbb P\,\(\,#1\,\)\,}

\newcommand{\cE}[2]{\mathbb{E}\,\(\,#1\,\big|\,#2\,\)\,}
\newcommand{\cP}[2]{\mathbb{P}\,\(\,#1\,\big|\,#2\,\)\,}

\newcommand{\qed}{\hfill$\square$}

\newcommand{\Proof}{\textsc{Proof. }}

\newcommand{\Rra}{R^{\text{random}}}
\newcommand{\Sra}{S^{\text{random}}}

\begin{document}

\title{Fourier-Walsh coefficients for a coalescing flow\\
  (discrete time)}
\author{Boris Tsirelson}
\date{}

\maketitle

\begin{abstract}
A two-dimensional array of independent random signs produces
coalescing random walks. The position of the walk, starting at the
origin, after $ n $ steps is a highly nonlinear, noise sensitive
function of the signs. A typical term of its Fourier-Walsh expansion
involves the product of $ \sim \sqrt n $ signs.
\end{abstract}

\section*{Introduction}

\begin{figure}[b]
\begin{center}
\setlength{\unitlength}{1cm}
\begin{picture}(12,2.6)
\put(-0.35,0.5){\includegraphics{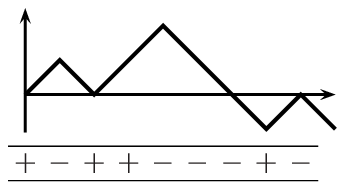}}
\put(1.5,0){\makebox(0,0){(a)}}
\put(6.0,0.2){\includegraphics{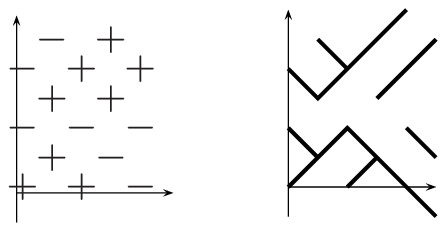}}
\put(8.4,0){\makebox(0,0){(b)}}
\end{picture}
\caption[]{\label{fig1}\small
(a) One-dimensional array of random signs produces a random walk. (b)
Two-dimensional array of random signs produces coalescing walks.}
\end{center}
\end{figure}

The simple random walk is driven by a one-dimensional array of
independent random signs (Fig.~\ref{fig1}a). The walk is a
\emph{stable} function of the signs, in the sense that changing at
random a small fraction of the signs results in a small change of the
walk.

\begin{figure}[htbp]
\begin{center}
\setlength{\unitlength}{1cm}
\begin{picture}(13,10)
\put(0,5.8){\includegraphics[width=6cm,height=4.5cm]{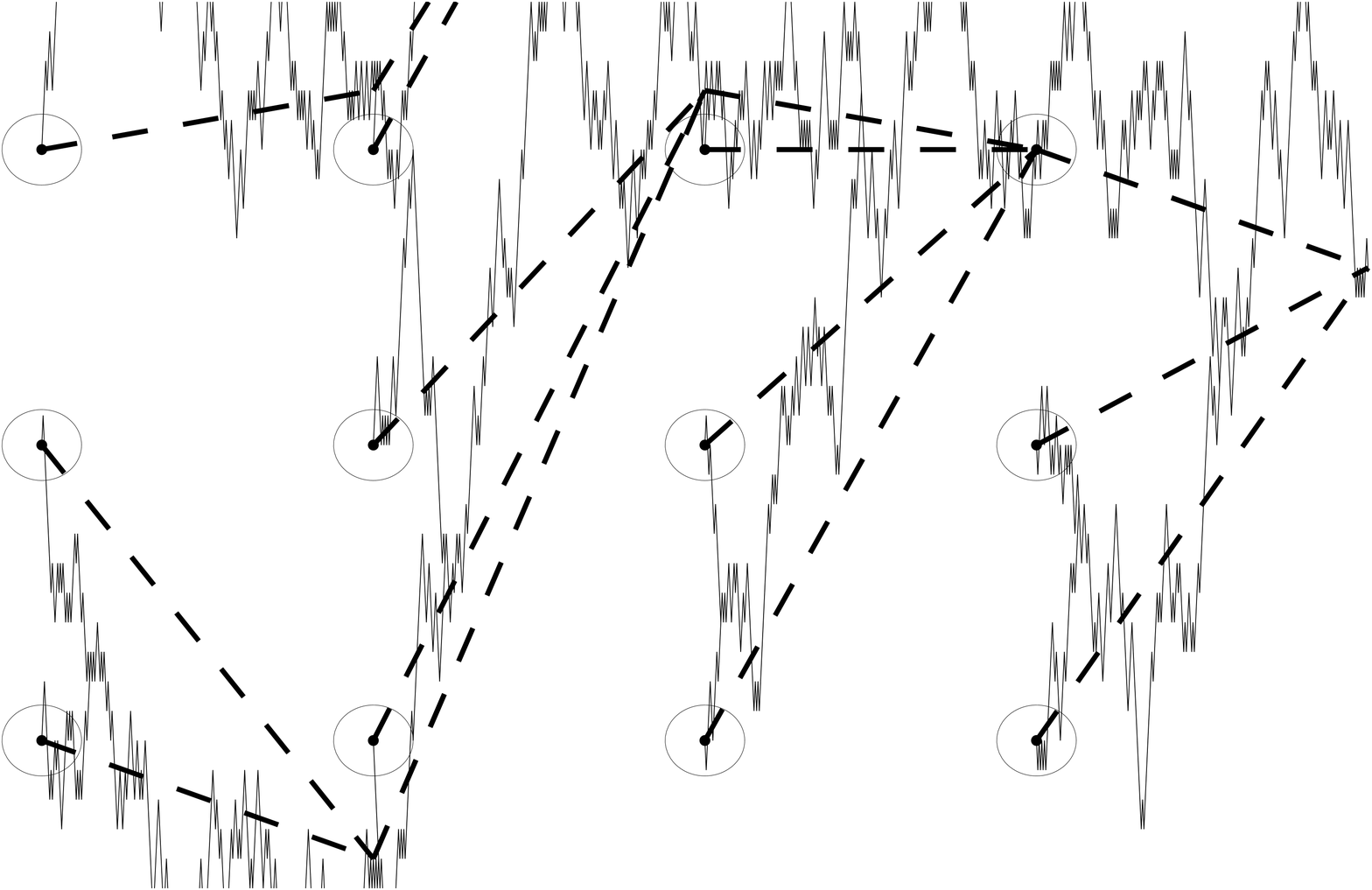}}
\put(7,5.8){\includegraphics[width=6cm,height=4.5cm]{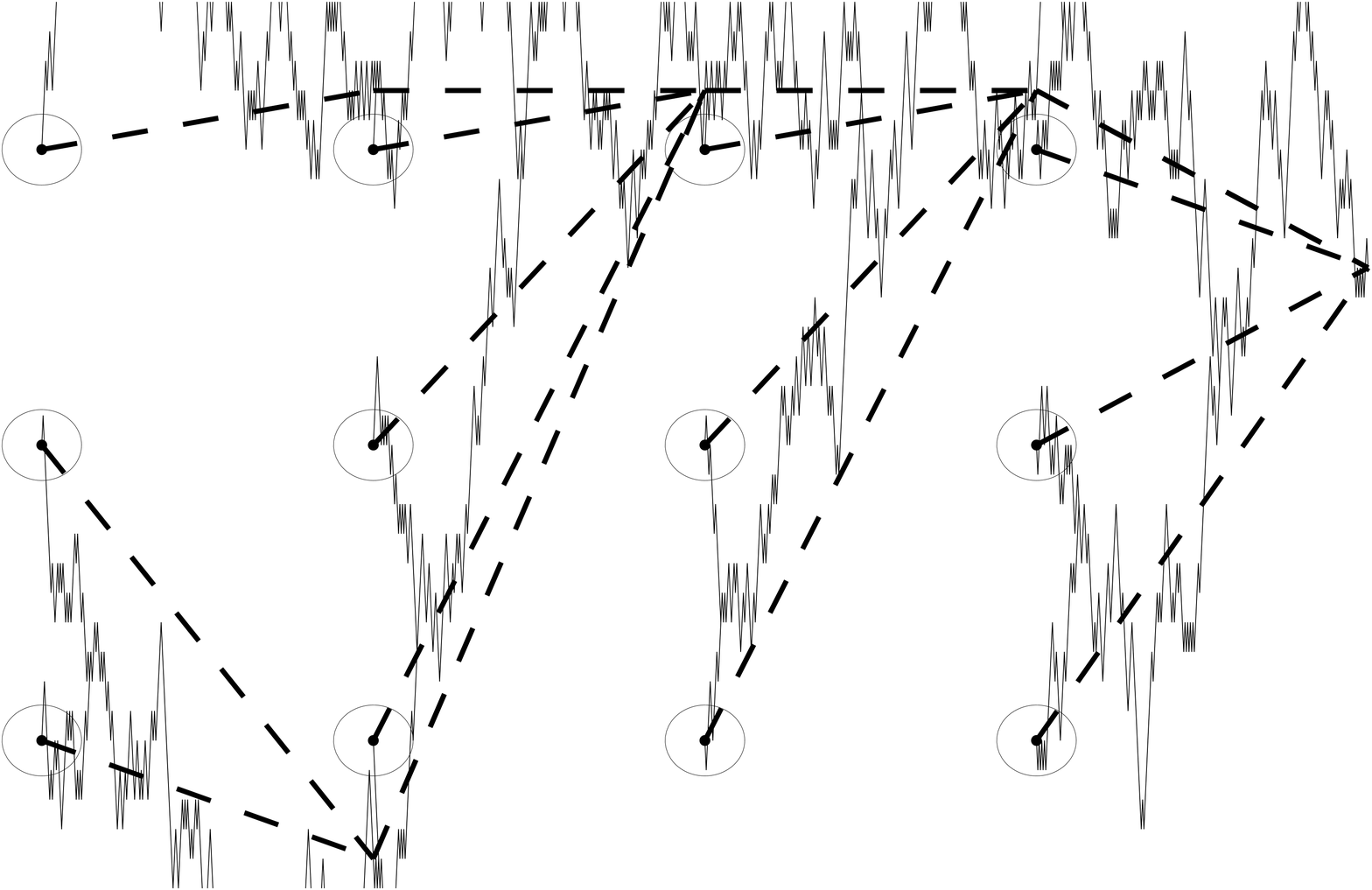}}
\put(0,0.3){\includegraphics[width=6cm,height=4.5cm]{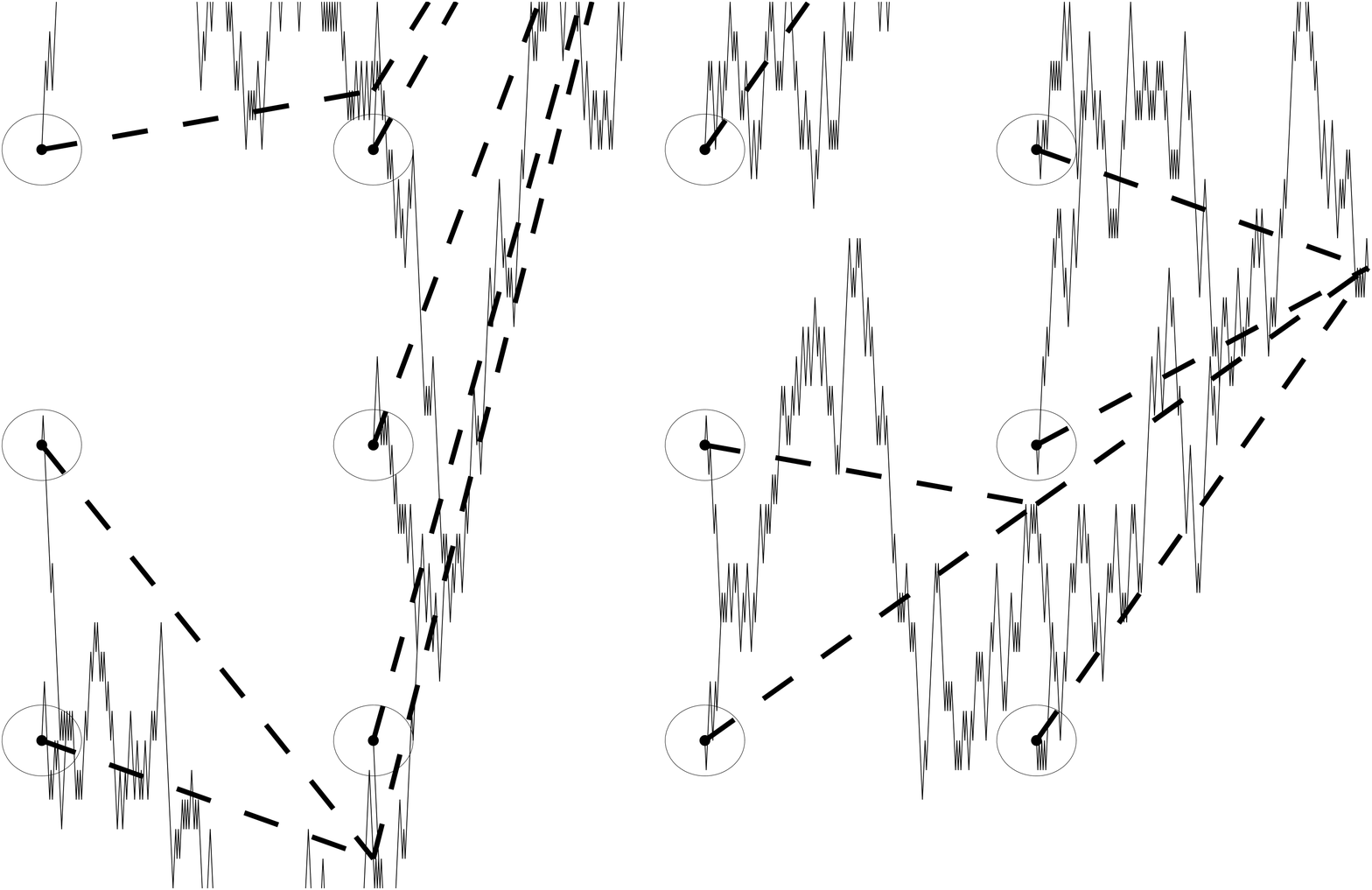}}
\put(7,0.3){\includegraphics[width=6cm,height=4.5cm]{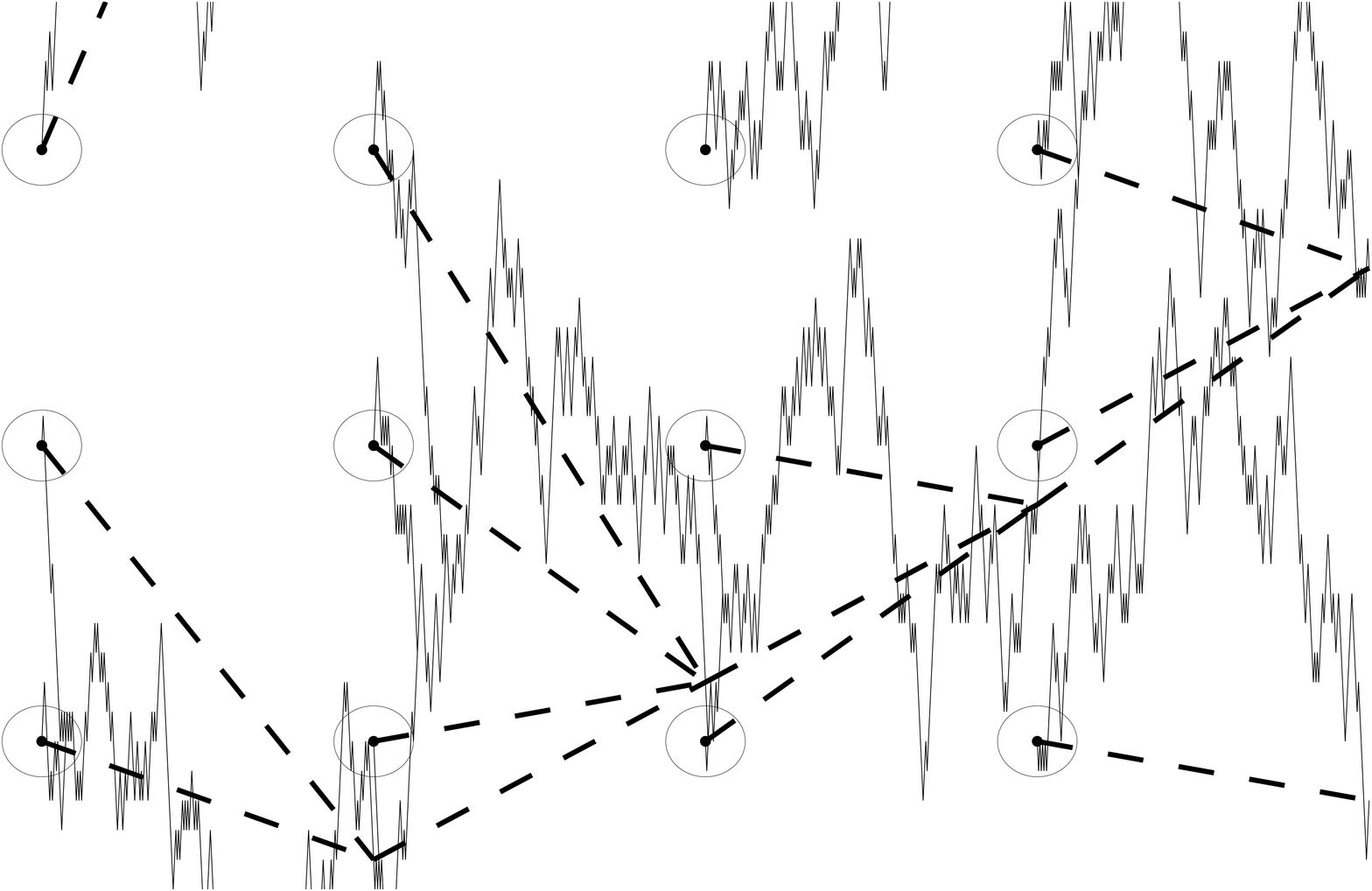}}
\put(3,5.5){\makebox(0,0){(a)}}
\put(10,5.5){\makebox(0,0){(b)}}
\put(3,0){\makebox(0,0){(c)}}
\put(10,0){\makebox(0,0){(d)}}
\end{picture}
\caption[]{\label{fig2}\small
Noise sensitivity. Coalescing walks on the grid $ 1000\times30 $,
starting from $ 4\times3 = 12 $ points (circled). Unperturbed array of
random signs (a). Perturbed array: each random sign is flipped with
probability 0.025 (b). Further perturbation of the same type
(c). Still further (d).}
\end{center}
\end{figure}

A two-dimensional array of independent random signs may be used for
producing the simple system of coalescing random walks
(Fig.~\ref{fig1}b). Changing at random a small fraction of the signs
causes a dramatic change of the walks (Fig.~\ref{fig2}). These are
\emph{sensitive} functions of the signs.

\begin{figure}[htbp]
\begin{center}
\setlength{\unitlength}{1cm}
\begin{picture}(10,5.5)
\put(-0.5,0.6){\includegraphics{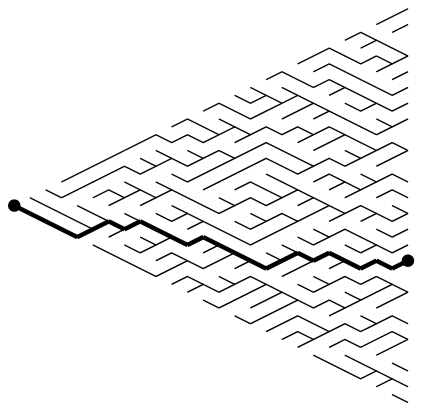}}
\put(2,0){\makebox(0,0){(a)}}
\put(5.75,0.4){\includegraphics{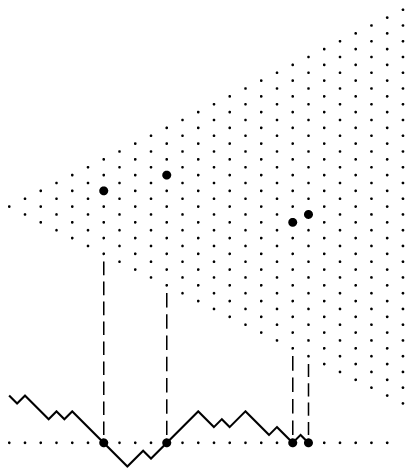}}
\put(8,0){\makebox(0,0){(b)}}
\end{picture}
\caption[]{\label{fig3}\small
The function of a triangular array of random signs (a), and
its typical spectral set (b).}
\end{center}
\end{figure}

\enlargethispage*{4pt}

Stability and sensitivity, introduced by Benjamini, Kalai, Schramm
\cite{BKS}, can be formulated also in terms of the Fourier-Walsh
transform. Any function of random signs can be written as a
polynomial, each term being a product of several signs, with a
coefficient. These are Fourier-Walsh coefficients. For a stable
function, they are concentrated on low degrees (frequencies), for a
sensitive function --- on high frequencies.

The simplest nontrivial function of coalescing walks is the final
position of the walk starting at the origin (Fig.~\ref{fig3}a).  Its
Fourier-Walsh transform is studied here. A typical term of the
polynomial appears to involve about $ \sqrt n $ signs (which is more
than enough for sensitivity). The corresponding set on the time axis
is distributed like the set of zeros of the simple random walk in
twice faster time (Fig.~\ref{fig3}b).  Implications for continuous
time will be published separately.

\vspace{-4pt}
\section{The model and results}\label{sect1}

\begin{sloppypar}
Choose $ n \in \{ 1,2,\dots \} $ and introduce a triangular array of
independent random signs
\[
\tau (x,y), \quad x = 0,1,\dots,n-1, \quad y = -x, -x+2, \dots, x-2,
x,
\]
each $ \tau (x,y) $ being $ -1 $ or $ +1 $ with probabilities $ 1/2
$. Define $ W (x) $ for $ x = 0,1,\dots,n $ by
\[
W(0) = 0; \quad W(x+1) = W(x) + \tau (x, W(x)) \, .
\]
Clearly, $ W(\cdot) $ is a simple random walk. In particular, $ \E \,
W(n) = 0 $ and $ \E \, W^2 (n) = n $. We consider the Fourier-Walsh
transform of the function $ \xi = n^{-1/2} W(n) $ of our random signs
$ \tau(x,y) $:
\[\begin{split}
& \frac1{\sqrt n} W(n) = \sum_{S\subset I} \hat\xi(S) \tau(S) \, , \\
& I = \{ (x,y) : x=0,1,\dots,n-1, \, y = -x, -x+2, \dots, x-2, x \} \,
 , \\
& \tau(S) = \prod_{(x,y)\in S} \tau(x,y) \, , \quad \hat\xi(S) = \Ex{
  \tau(S) \xi } = \frac1{\sqrt n} \Ex{ \tau(S) W(n) } \, .
\end{split}\]
The equality $ \E \, \xi = 0 $ becomes $ \hat\xi (\emptyset) = 0 $,
while $ \E \, \xi^2 = 1 $ becomes $ \sum_{S\subset I} | \hat\xi(S) |^2
= 1 $, since $ \tau(S) $ are an orthonormal basis of $ l_2 \(
\{-1,+1\}^I \) $.
\end{sloppypar}

\noindent\dotfill\vspace{-6pt}

\begin{proposition}\label{prop1.1}
$ \hat\xi (S) = 0 $ unless $ S $ is of the form
\begin{equation}\label{1.2}\begin{split}
& S = \{ (x_1,y_1), \dots, (x_m,y_m) \} \, , \\
& m \in \{ 1, \dots, n \} \, , \\
& 0 \le x_1 < \dots < x_m < n \, , \\
& | y_1 | \le x_1 \, , \\
& | y_{k+1} - y_k | \le x_{k+1} - x_k \quad \text{for } k=1,\dots,m-1
  \, ,
\end{split}\end{equation}
in which case
\[
\hat\xi(S) = \frac1{\sqrt n} p(x_1,y_1) q(S) \, ,
\]
where
\[\begin{split}
& q(S) = \\
& = \prod_{k=1}^{m-1} \frac{ p (x_{k+1}-x_k-1,
  y_{k+1}-y_k-1) - p (x_{k+1}-x_k-1, y_{k+1}-y_k+1) }{ 2 } \, , \\
& p(x,y) = 2^{-x} \frac{ x! }{ \left( \frac{x-y}2
  \right) ! \left( \frac{x+y}2 \right) ! } \, .
\end{split}\]
\noindent\dotfill
\end{proposition}

Define a probability distribution $ \mu_\xi $ (the spectral measure of
$ \xi $) on the set $ 2^I = \{ S : S \subset I \} $ by
\[
\mu_\xi (A) = \sum_{S\in A} | \hat\xi(S) |^2 \quad \text{for } A
\subset 2^I \, .
\]
Let $ \Sra $ be a random subset of $ I $, distributed according to $
\mu_\xi $. We know that $ \Sra $ is of the form \eqref{1.2} with
probability $ 1 $. The projection of $ \Sra $ onto the first axis is
another random set $ \Rra = \{ x_1,\dots,x_m \} $. The following
result describes the probability distribution of $ \Rra $.

\noindent\dotfill\vspace{-6pt}

\begin{proposition}\label{prop1.2}
Whenever $ m \in \{1,\dots,n\} $ and $ 0 \le x_1 < \dots < x_m < n $,
\begin{multline*}
\Pro{ \Rra = \{ x_1, \dots, x_m \} } = \\
\quad = \frac1n p(2x_1,0) \cdot \prod_{k=1}^{m-1} \(
  p(2x_{k+1}-2x_k-2,0) - p (2x_{k+1}-2x_k,0) \) \, .
\end{multline*}
\noindent\dotfill
\end{proposition}

It means that $ \Rra $ may be described as follows. First, we choose
its maximal element $ x_m $ (without knowing $ m $ for now) uniformly
on $ \{ 0,\dots,{n-1} \} $. Then we introduce a simple random walk $ V =
\( V(0), V(1), \dots \) $ and consider $ V (2(x-x_m)) $; that may be
thought of as a random walk starting at $ (x_m,0) $ and going backward
in time, twice faster than usual, see Fig.~\ref{fig3}b. Its zeros are
just $ \Rra $,
\begin{equation}\label{1.4}
\Rra = \{ x \in \{0,\dots,x_m\} : V (2(x-x_m)) = 0 \} \, .
\end{equation}
It is strange! After proving Prop.~\ref{prop1.2} we may introduce $ V
$ satisfying \eqref{1.4}. It would be more natural to do other way
round.

\begin{problem}
Can Prop.~\ref{prop1.2} be deduced from \eqref{1.4}? That is, can $ V
$ be introduced somehow before proving Prop.~\ref{prop1.2}?
\end{problem}

\section{Proofs}\label{sect2}

The random path $ W = \( W(0), \dots, W(n) \) $ determines a partition
of our probability space $ \Om = \{ -1, +1 \}^I $ into subsets $ \{ W
= w \} $ indexed by paths $ w $, that is, by sequences $ w = \( w(0),
\dots, w(n) \) $ such that $ w(0) = 0 $ and $ w(k+1) - w(k) = \pm1 $
for $ k = 0, \dots, n-1 $.

\sk

\textsc{Proof of Prop.~\ref{prop1.1}.}\\*
We have $ \hat\xi(S) = \Ex{ \tau(S) \xi } = \sum_w \Ex{ \tau(S) \xi
\mathbf1_{(W=w)} } $, where $ \mathbf1_{(W=w)} $ is the indicator of
the event $ W = w $. However, $ W(n) \mathbf1_{(W=w)} = w(n)
\mathbf1_{(W=w)} $, thus $ \hat\xi(S) = n^{-1/2} \sum_w w(n) \Ex{
\tau(S) \mathbf1_{(W=w)} } $. Clearly, $ \mathbf1_{(W=w)} $ is
independent of all $ \tau(x,y) $ with $ w(x) \ne y $. Therefore $ \Ex
{ \tau(S) \mathbf1_{(W=w)} } = 0 $ unless $ w $ passes through all
points of $ S $, which proves \eqref{1.2}.

\begin{sloppypar}
Now, $ S = \{ (x_1,y_1), \dots, (x_m,y_m) \} $, $ 0 \le x_1 < \dots <
x_m < n $. Let $ w $ pass through all points of $ S $, then $
\mathbf1_{(W=w)} \tau(S) = \mathbf1_{(W=w)} \prod_{k=1}^m \( { w(x_k+1)
- y_k } \) $ and $ \Ex{ \tau(S) \mathbf1_{(W=w)} } = 2^{-n}
\prod_{k=1}^m \( w(x_k+1) - y_k \) $. Therefore $ \hat\xi(S) =
n^{-1/2} 2^{-n} \sum_w w(n) \prod \( w(x_k+1) - y_k \) $, the sum
being taken over $ w $ that satisfy $ w(x_k) = y_k $ for all $ k $. 
The set of such $ w $ is a product of $ m+1 $ sets, corresponding to
the partition of $ [0,n) $ into $ [0,x_1) $, $ [x_1,x_2) $, \ldots,
$ [x_{m-1},x_m) $, $ [x_m,n) $, and the sum over $ w $ decomposes into
a product of $ m+1 $ sums. The first sum is equal to $ \binom{ x_1 }{
(x_1-y_1)/2 } = 2^{x_1} p(x_1,y_1) $. The second sum, containing $
w(x_1+1) - y_1 $, is equal to $ 2^{x_2-(x_1+1)} p \( x_2 - (x_1+1),
y_2 - (y_1+1) \) - 2^{x_2-(x_1+1)} p \( x_2 - { (x_1+1) }, y_2 - (y_1-1)
\) $. And so on. The last sum, containing $ w (x_m+1) - y_m $, but
also $ w(n) $, is equal to $ 2^{n-(x_m+1)} (y_m+1) - 2^{n-(x_m+1)}
(y_m-1) = 2^{n-x_m} $. So,
\begin{multline*}
\hat\xi(S) = \frac1{\sqrt n} \cdot 2^{-n} \cdot 2^{x_1} p(x_1,y_1)
  \times \\
\bigg( \prod_{k=1}^{m-1} 2^{x_{k+1}-x_k-1} \(
  p(x_{k+1}-x_k-1,y_{k+1}-y_k-1) - p( \dots, y_{k+1}-y_k+1) \)
  \bigg)  \\
\times 2^{n-x_m} = \frac1{\sqrt n} p(x_1,y_1) q(S) \, .
\end{multline*}
\qed
\end{sloppypar}

\sk

The following general fact holds for arbitrary $ I $ and $ \xi $.

\begin{lemma}\label{lm1}
Let $ I $ be a finite set, $ \( \tau_i \)_{i\in I} $ be independent
random signs ($ \E \tau_i = 0 $), and $ \xi = \sum_{S\subset I}
\hat\xi(S) \tau(S) $, where $ \tau(S) = \prod_{i\in S} \tau_i $. Then
for any $ E \subset I $ and $ T \subset I \setminus E $,
\[
\sum_{S\subset E} \big| \hat\xi (S \cup T) \big|^2 = \E \, \big| \,
\cE{ \tau(T)\xi }{ \tau|_E } \! \big|^2
\]
(the conditional expectation $ \cE{\dots}{\tau|_E} $ is taken
w.r.t.\ all $ \tau(i) $ for $ i \in E $).
\end{lemma}

\textsc{Proof.}
We have
\[
\xi = \sum_{S_1\subset E} \sum_{S_2\subset I\setminus E} \hat\xi
(S_1\cup S_2) \tau(S_1) \tau(S_2) \, ,
\]
therefore
\[
\cE{ \tau(T)\xi }{ \tau|_E } = \sum_{S_1\subset E} \tau(S_1)
  \sum_{S_2\subset I\setminus E} \hat\xi (S_1\cup S_2)
  \, \cE{ \tau(T) \tau(S_2) }{ \tau|_E } \, .
\]
However, $ \cE{ \tau(T) \tau(S_2) }{ \tau|_E } = \Ex{ \tau(T)
\tau(S_2) } = 0 $ unless $ S_2 = T $. So, $ \cE{ \tau(T)\xi }{
\tau|_E } = \sum_{S_1\subset E} \tau(S_1) \hat\xi ( S_1 \cup T ) $ and
$ \E \, | \dots |^2 = \sum_{S_1\subset E} | \hat\xi ( S_1 \cup T ) |^2
$.\qed

\medskip

We return to special $ I $ and $ \xi $ introduced in
Sect.~\ref{sect1}. Recall the projection $ \Rra $ of the random set $
\Sra $. It is never empty. The following result describes the
(uniform, in fact) distribution of the maximal element of $ \Rra $.

\begin{lemma}\label{lm2}
$ \Pro{ \Rra \subset [0,k] } = \dfrac{k+1}n $ for $ k = 0,\dots,n-1 $.
\end{lemma}

\textsc{Proof.}
Lemma \ref{lm1} for $ E = \{ (x,y) \in I : x \le k \} $ and $ T = \emptyset $
gives $ \Pro{ \Rra \subset [0,k] } \! = \E \, \big| \, \cE{
\xi }{ \tau|_E } \big|^2 $. However, $ \cE{ W(n) }{ \tau|_E } \! = { W(k+1) }
$, thus, $ \Pro{ \Rra \subset [0,k] } = \E \, | n^{-1/2}
W(k+1) |^2 = \frac{k+1}n $.\qed

\medskip

Recall $ q(S) $ defined in Prop.~\ref{prop1.1} for $ S $ of the form
\eqref{1.2}.

\begin{lemma}\label{lm3}
Let $ k \in \{ 0,\dots,n-2 \} $, $ E = \{ (x,y) \in I : x \le k \} $,
and $ S \subset I \setminus E $ be of the form \eqref{1.2}, then
\[
\cE{ \tau(S)\xi }{ \tau|_E } = \frac1{\sqrt n} p \( x_1 - (k+1), y_1 -
W(k+1) \) q(S)
\]
($ x_1, y_1 $ being defined by \eqref{1.2}).
\end{lemma}

\Proof
Similarly to the proof of Prop.~\ref{prop1.1}, but now $ w(0),
\dots, w(k+1) $ are fixed, and the sum is taken over $ w(k+2), \dots,
w(n) $.

\medskip

\begin{lemma}\label{lm4}
Let $ k $, $ E $, $ S $ be as in Lemma \ref{lm3}, then
\[
\Pro{ \Sra \setminus E = S } = \frac1n \bigg( \sum_y p(k+1,y) p^2
(x_1-k-1, y_1-y) \bigg) q^2(S) \, .
\]
\end{lemma}

\textsc{Proof.} Lemma \ref{lm1} gives
\[
\Pro{ \Sra \setminus E = S } = \E \, \big| \, \cE{ \tau(S)\xi }{
\tau|_E } \big|^2 \, .
\]
Lemma \ref{lm3} gives
\[
\cE{ \tau(S)\xi }{ \tau|_E } = \frac1{\sqrt n} p ( x_1 - (k+1), y_1 -
W(k+1) ) q(S) \, .
\]
So, $ \Pro{ \Sra \setminus E = S } = \E \, |\dots|^2 = \frac1n q^2(S)
\, \E \, p^2 ( x_1 - (k+1), y_1 - W(k+1) ) $. It remains to note that
$ \E \, p^2 ( x_1 - (k+1), y_1 - W(k+1) ) = \sum_y p^2 ( x_1-k-1,
y_1-y ) p (k+1,y) $.\qed

\medskip

Given a set $ S \subset I $ of the form \eqref{1.2}, we introduce the
set $ S^\updownarrow $ of all vertical shifts of $ S $. That is, for $
S = \{ (x_1,y_1), \dots, (x_m,y_m) \} $, $ S^\updownarrow $ consists
of the sets $ \{ (x_1,y_1+\de), \dots, (x_m,y_m+\de) \} $, where $ \de
$ satisfies $ y_1+\de \in \{ -x_1, -x_1+2, \dots, x_1-2, x_1 \} $.

\begin{lemma}\label{lm5}
Let $ k $, $ E $, $ S $ be as in Lemma \ref{lm3}, then
\[
\Pro{ \Sra \setminus E \in S^\updownarrow } = \frac1n p \(
2(x_1-k-1), 0 \) q^2(S) \, .
\]
\end{lemma}

\Proof
By Lemma \ref{lm4}, the probability is
\[
\sum_{y_1} \frac1n \bigg( \sum_y p(k+1,y) p^2 (x_1-k-1, y_1-y) \bigg) q^2(S)
\, ,
\]
since $ q(\cdot) $ is shift-invariant. However,
\begin{multline*}
\sum_{y,y_1} p(k+1,y) p^2 (x_1-k-1, y_1-y) = \\
= \bigg( \sum_y p(k+1,y) \bigg) \cdot \bigg(
  \sum_y p^2 (x_1-k-1, y) \bigg) = \\
= \sum_y p^2 (x_1-k-1, y) = p \( 2(x_1-k-1), 0 \) \, ;
\end{multline*}
the latter equality is a well-known property of binomial coefficients, see
for instance \cite[II.12.11]{Feller}.\qed

\begin{lemma}\label{lm6}
Let $ 0 \le k < x_1 < \dots < x_m < n $, then
\begin{multline*}
\cP{ \Rra \cap [k,x_1) = \emptyset }{ \Rra \cap [x_1,n) =
  \{x_1,\dots,x_m\} } = \\
= p ( 2(x_1-k), 0 ) \, .
\end{multline*}
\end{lemma}

\textsc{Proof.}
The event $ \Rra \cap [x_1,n) = \{x_1,\dots,x_m\} $ is the union of
disjoint events of the form $ \Sra \setminus E_1 \in S^\updownarrow
$, where $ E_1 = \{ (x,y) \in I : x<x_1 \} $ and $ S \subset I
\setminus E_1 $, $ S = \{ (x_1,y_1), \dots, (x_m,y_m) \} $ for some $
y_1, \dots, y_m $ such that $ S $ satisfies \eqref{1.2}. Given such $
S $, we use Lemma \ref{lm5}, its $ k $ being our $ x_1 - 1 $; we get
\[
\Pro{ \Sra \setminus E_1 \in S^\updownarrow } = \frac1n q^2 (S)
\, ,
\]
since $ p(0,0) = 1 $.

We use Lemma \ref{lm5} once again, its $ k $ being now equal to our $
k-1 $; we get
\[
\Pro{ \Sra \setminus E_2 \in S^\updownarrow } = \frac1n p \(
2(x_1-k), 0 \) q^2 (S) \, ,
\]
where $ E_2 = \{ (x,y) \in I : x < k \} $.

Note that $ \Sra \setminus E_2 \in S^\updownarrow $ if and only if $
\Sra \setminus E_1 \in S^\updownarrow $ and $ \Rra \cap [k,x_1) =
\emptyset $. Thus,
\begin{multline*}
\cP{ \Rra \cap [k,x_1) = \emptyset }{ \Sra \setminus E_1 \in
  S^\updownarrow } = \\
= \frac{ \Pro{ \Sra \setminus E_2 \in S^\updownarrow } }{ \Pro{
   \Sra \setminus E_1 \in S^\updownarrow } } = p \( 2(x_1-k), 0 \)
  \, .
\end{multline*}
The conditional probability does not depend on $ S $. Summing over all
$ S $ (one $ S $ within each equivalence class $ S^\updownarrow $), we
get
\begin{multline*}
\cP{ \Rra \cap [k,x_1) = \emptyset }{ \Rra \cap [x_1,n) = \{
 x_1,\dots,x_m\} } = \\
= p \( 2(x_1-k), 0 \) \, .
\end{multline*}
\qed

\sk

\textsc{Proof of Prop.~\ref{prop1.2}.}
We have
\begin{multline*}
\Pro{ \Rra = \{ x_1,\dots,x_m \} } = \\
= \Pro{ \Rra \cap [x_m,n) = \{ x_m \} } \times \\
\times \cP{ \Rra \cap [x_{m-1}, x_m) = \{x_{m-1}\} }{ \Rra \cap
  [x_m,n) = \{x_m\} } \cdot \ldots \times \\
\times \cP{ \Rra \cap [x_1, x_2) = \{x_1\} }{ \Rra \cap
  [x_2,n) = \{x_2,\dots,x_m\} } \times \\
\times \cP{ \Rra \cap [0, x_1) = \emptyset }{ \Rra \cap [x_1,n) =
  \{x_1,\dots,x_m\} } \, .
\end{multline*}
The first factor is
\begin{multline*}
\Pro{ \Rra \cap [x_m,n) = \{ x_m \} } = \\
= \Pro{ \Rra \subset [0,x_m] } - \Pro{ \Rra \subset [0,x_m)
  } = \\
= \frac{x_m+1}n - \frac{x_m}n = \frac1n
\end{multline*}
by Lemma \ref{lm2}. The second factor is
\begin{multline*}
\cP{ \Rra \cap [x_{m-1}, x_m) = \{x_{m-1}\} }{ \Rra \cap
  [x_m,n) = \{x_m\} } = \\
\cP{ \Rra \cap (x_{m-1}, x_m) = \emptyset }{ \dots } - \cP{ \Rra
  \cap [x_{m-1}, x_m) = \emptyset }{ \dots } \\
= p \( 2(x_m-x_{m-1}-1), 0 \) - p \( 2(x_m-x_{m-1}), 0 \)
\end{multline*}
by Lemma \ref{lm6}. And so on. The last factor is
\[\begin{split}
& \cP{ \Rra \cap [0, x_1) = \emptyset }{ \Rra \cap [x_1,n) =
  \{x_1,\dots,x_m\} } = p ( 2x_1, 0 ) \, .
\end{split}\]
So,
\begin{multline*}
\Pro{ \Rra = \{x_1,\dots,x_m\} } = \frac1n \times \\
\times \bigg( \prod_{k=1}^{m-1} \( p ( 2(x_{k+1}-x_k-1), 0 ) - p (
  2(x_{k+1}-x_k), 0 ) \) \bigg) \times \\
\times p ( 2x_1, 0 ) \, .
\end{multline*}
\qed

\bigskip
\filbreak
\begingroup
{
\small
\begin{sc}
\parindent=0pt\baselineskip=12pt
\def\emailwww#1#2{\par\qquad {\tt #1}\par\qquad {\tt #2}\medskip}

School of Mathematics, Tel Aviv Univ., Tel Aviv
69978, Israel
\emailwww{tsirel@math.tau.ac.il}
{http://math.tau.ac.il/$\sim$tsirel/}
\end{sc}
}
\filbreak

\endgroup

\end{document}